\documentclass[12pt]{article}
\usepackage{amsmath,amssymb,amsthm}
\usepackage[russian]{babel}
  \def\R{{\mathbb R}} 
\def\C{{\mathbb C}}

\long\def\comment#1\endcomment{}

\begin{document}

\smallskip
\centerline{\bf On some results by S. Belkov and I. Korepanov}
\smallskip
\centerline{\bf A. Skopenkov}

\bigskip
This note is purely expositional and is a complement to math review MR2730150 to [BK].
%the reviewed paper being correct or new.
I tried to present in a clear way the statement of Theorem 1 and main definition from [BK].
I bear no responsibility for results of [BK].

In [BK] the authors consider purely mathematical problems of finding an invariant of a
3-manifold (Theorem 4; Theorem 2 is a lemma for Theorem 4) and finding certain identities
(Theorem 1 and its matrix version Theorem 3).
No physical consequences are presented (this is worth mentioning because
the paper is published in the journal whose title is translated as `Theoretical and
mathematical physics').
Main definitions are not clearly presented.
\footnote{It is not written whether we construct an invariant of a pair $(M,C)$ of a 3-manifold
and $C$, or we construct an invariant of $M$ using $C$ and then prove the independence of $C$.
It is not clear whether it is constructed a family of invariants depending on $n$, or only one
invariant (whose independence of $n$ is proved).
Definition of $f_2$ in p. 518 is meaningless because $\zeta_i$ are not defined in \S3.
Presumably in all terms of (5) except the left $\C^{nN_0'}$, $\C$ should be changed to the field
$\C(\zeta_1,\dots,\zeta_{N_0})$ of rational functions; it is not clear whether the left $\C^{nN_0'}$
should be changed like that.}
%If a reader would be enough interested to reconstruct the statements and the proofs, he/she can judge whether
%his/her interpretation of the results is interesting, new and well-known.
%In this note I present a reconstruction of Theorem 1 and a
%possible reconstruction of definition of $f_2,f_3$ and $f_4$; this definition is required
%for all other main results of the paper.

\smallskip
{\it A reconstruction of Theorem 1.}
A Grassmanian algebra (over $\R$ or $\C$) is an associative anti\-com\-mu\-ta\-ti\-ve algebra with
unity.
%, generators $a_i$ and relations $a_ia_j$
Let $G$ be a Grassmanian algebra with generators $\{a_s\}$.
Each element of $G$ can be represented as a polynomial of $\{a_s\}$ of degree at most 1 over each $a_s$.
The {\it Berezin integral} corresponding to $a_s$ is a linear operator $B_s:G\to G$ defined by $B_s(Pa_s+Q):=P$,
where $P$ and $Q$ are polynomials (of degree at most 1) of all the generators except $a_s$.
Define $B_s:G[\zeta_1,\dots,\zeta_5]\to G[\zeta_1,\dots,\zeta_5]$ by applying $B_s$ to all the coefficients
of a polynomial of $\zeta_1,\dots,\zeta_5$.
Consider Grassmanian algebra having generators $\{a_{ijk}\}$ numbered by 3-element subsets $\{i,j,k\}\subset\{1,2,3,4,5\}$
(and, possibly, other generators).
For a 4-element subset $\{i,j,k,l\}\subset\{1,2,3,4,5\}$ define
$$f_{ijkl}:=\sum\limits_{\{p,q\}\subset \{i,j,k,l\}}(\zeta_p-\zeta_q)a_{pqr}a_{pqs}\in G[\zeta_1,\dots,\zeta_5],$$
%\zeta_{ij}a_{ijk}a_{ijl}-\zeta_{ik}a_{ikj}a_{ikl}+\zeta_{ik}a_{ikj}a_{ikl}+\zeta_{jk}a_{jki}a_{jkl}-\zeta_{jl}a_{jli}a_{jl}\zeta_{ij}a_{ijk}a_{ijl}$
where in the summand we take any {\it ordered} pair $(p,q)$ corresponding to $\{p,q\}$ and define $r,s$ so that
$\{p,q,r,s\}=\{i,j,k,l\}$ and $p,q,r,s$ is an {\it even} permutation of $i,j,k,l$.
Theorem 1 states that
$$B_{123}(f_{1234}f_{1235})=B_{145}B_{245}B_{345}(f_{1245}f_{2345}f_{1345}).$$
\quad
{\it A possible reconstruction of definition of $f_2,f_3$ and $f_4$ (this definition is required
for all other main results of [BK]).}
Take a triangulation of a compact 3-manifold.
Denote by $N_0', N_2'$ and $N_3$ the number of interior vertices, of interior 2-simplices and of all 3-simplices, respectively.
Since $N_2'\le2N_3$, we can take a subset $X$ of the set of boundary 2-simplices such that
$\#X=2N_3-N_2'$.
Let $G:=\C(\zeta_1,\dots,\zeta_{N_0'})$ be the field of rational functions.
Let us define maps from the diagram
$$G^{nN_0'}\overset{f_2}\to G^{2nN_3}\overset{f_3}\to G^{2nN_3}\overset{f_4}\to G^{nN_0'}.$$
Represent an element of $G^{nN_0'}$ as $(u_1,\dots,u_{N_0'})$, $u_s\in G^n$.
Take an ordering on the union $Y$ of $X$ with the set of interior 2-simplices.
Represent an element of the left $G^{2nN_3}$ as a vector $\vec v$ with components $v_{ijk}\in G^n$, $\{i,j,k\}\in Y$ and $i<j<k$.
%with $\C^n$-components indexed by 2-simplices of $Y$.
Define $f_2$ by
$$[f_2(v_1,\dots,v_{N_0'})]_{ijk}:=\frac{v_i-v_j}{\zeta_i-\zeta_j}-\frac{v_i-v_k}{\zeta_i-\zeta_k}.$$
Take an ordering on the set of all 3-simplices.
Represent an element of the right $G^{2nN_3}$ as a vector $\vec w$ with components $w_{i,ijkl},w_{j,ijkl}\in G^n$, $\{i,j,k,l\}$ is a 3-simplex, $i<j<k<l$.
Define $f_3$ by
$$[f_3(\vec v)]_{p,ijkl}:=
\begin{cases}v_{ijk}+v_{ikl}-v_{ijl}, &p=i\\
v_{ijk}\frac{\zeta_i-\zeta_k}{\zeta_k-\zeta_j}-v_{ijl}\frac{\zeta_i-\zeta_l}{\zeta_l-\zeta_j}-v_{jkl}, &p=j\end{cases}.$$
Given $w_{i,ijkl},w_{j,ijkl}\in G^n$ define $w_{k,ijkl},w_{l,ijkl}\in G^n$ by the system of equations
$$w_{i,ijkl}+w_{j,ijkl}+w_{k,ijkl}+w_{l,ijkl}=0,\quad \zeta_iw_{i,ijkl}+\zeta_jw_{j,ijkl}+\zeta_kw_{k,ijkl}+\zeta_lw_{l,ijkl}=0.$$
Define $f_4$ by
$$[f_4(\vec w)]_i:=\sum\limits_{i\text{ is a vertex of 3-simplex }(p,q,r,s),\ p<q<r<s}\varepsilon(p,q,r,s)w_{i,pqrs},$$
where  $\varepsilon(p,q,r,s)$ is 1 or $-1$ according to the orientation of 3-simplex $(p,q,r,s)$ positive or not.
Theorem 2 states that $f_3\circ f_2=0$ and $f_4\circ f_3=0$.

\smallskip
[BK] Bel'kov, S. I.; Korepanov, I. G. Matrix solution of the pentagon equation with
anti\-com\-mu\-ting
variables, Teoret. i Matemat. Fizika, 163:3 (2010), 513--528.

\end{document}